\documentclass{svproc}
\usepackage{epsfig,amsmath,amsfonts}
\usepackage{url}

\newcommand{\dd}{{\rm\,d}}
\DeclareMathOperator{\tr}{tr}
\newcommand{\eps}{\varepsilon}
\newcommand{\mymat}[1]{\mathbf{#1}}

\begin{document}
\mainmatter              % start of a contribution
\title{Stability boundary approximation\newline of periodic dynamics}
\titlerunning{Stability boundary approximation of periodic dynamics}  % abbreviated title (for running head)
\author{Anton O. Belyakov\inst{1}  %
\and Alexander P. Seyranian\inst{2}}
\authorrunning{Anton O. Belyakov and Alexander P. Seyranian} % abbreviated author list (for running head)
\tocauthor{Anton O. Belyakov and Alexander P. Seyranian}
\institute{Moscow School of Economics, Lomonosov Moscow State  University,\\
Leninskie Gory 1-61, 119234 Moscow, Russia;\\
National Research Nuclear University ``MEPhI'', Moscow, Russia;\\
Central Economics and Mathematics Institute, Russian Academy of Sciences;\\
\email{belyakov@mse-msu.ru}
%\email{I.Ekeland@princeton.edu},\\ WWW home page: \texttt{http://users/\homedir iekeland/web/welcome.html}
\and
Institute of Mechanics, Lomonosov Moscow State  University,\\
Michurinskiy prospect 1, 119192  Moscow, Russia,\\
\email{seyran@imec.msu.ru}}

\maketitle              % typeset the title of the contribution

\begin{abstract}
We develop here the method for obtaining approximate stability boundaries in the space of parameters for systems with parametric excitation. The monodromy (Floquet) matrix of linearized system is found  by averaging method. For system with 2 degrees of freedom (DOF) we derive general approximate stability conditions.  We study domains of stability with the use of fourth order approximations of monodromy matrix on example of inverted position of a pendulum with vertically oscillating pivot. Addition of small damping shifts the stability boundaries upwards, thus resulting to both stabilization and destabilization effects.
\end{abstract}
\section{Introduction}
Let us study the stability of the equilibrium $y\equiv 0$ of a nonlinear system governed by ordinary differential equations
$\dot y = \mathbf J(t)\, y + f(t,y)$, where for nonlinear vector-function $f$  exist constants $c$ and $\nu > 1$, such that $|f(t,y)|\leq c\, |y|^{\nu}$ for all $t\geq 0$ and $y\in \mathbb R^n$. Vector-function $f(t,y)$ and matrix $\mathbf J(t)$ are piecewise continuous in $t$. Moreover $\mathbf J(t)$ is bounded and $T$-periodic.
 
According to Lyapunov's theorem the trivial solution of such a nonlinear system is \emph{asymptotically stable} if all \emph{Lyapunov exponents} of the linear system $\dot{x}(t)=\mathbf J(t)\, x(t)$ are strictly negative and the solution is unstable if at least one Lyapunov exponent is strictly positive.\footnote{\emph{Lyapunov regularity} condition of the linear system holds when its matrix is periodic.} Asymptotic stability of the periodic linear system determines the asymptotic stability of the nonlinear system equilibrium and \emph{vice versa}.\footnote{Periodic linear system is asymptotically stable if, and only if, it is exponentially stable. Exponential stability of $\dot{x}(t)=\mathbf J(t)\, x(t)$ results in exponential stability, and hence asymptotic stability, of the nonlinear system solution $y\equiv 0$.}  Stability of  linear systems with time-periodic coefficients was also studied by Gaston Floquet \cite{Floquet}. Lyapunov exponent $\lambda$ of the linear periodic system can be expressed via its corresponding \emph{Floquet multiplier} $\rho$ as $\lambda = \frac{1}{T}\ln|\rho|$.  The theorem can also  be reformulated to compare absolute values of Floquet multipliers with 1. Floquet multipliers are the eigenvalues of the \emph{monodromy matrix} which is the fundamental matrix of the linear system taken at time $T$.

Monodromy matrix and Floquet multipliers can always be numerically calculated. Stability can also be checked by studying solutions of the system, see, e.g. \cite{XuWiercigroch2007,Butikov2017,BogMit} and references therein. But in practice it is often useful to have  analytical approximations of stability regions in parameter space. To obtain a straightforward technique for deriving such analytical stability boundary approximations of any order we combine Floquet theory with asymptotic method of \emph{averaging}, \cite{BogMit}.  

This technique yields same results as expansion of monodromy matrix in series in \cite{SeyMai} up to the terms of higher order than the order of approximation. These terms are automatically eliminated in the averaging scheme making the technique more convenient in practice.

We demonstrate the proposed technique on the example of inverted pendulum, where destabilizing effect of damping (shift of the lower stability boundary) is manifested in the fourth approximation, see \cite{Seyran}, though it would be natural to expect stabilization by damping, see, e.g. \cite{XuWiercigroch2007,ArkhipovaLuongo2016}. But there is also destabilizing effect  (shift of the upper boundary) demonstrated numerically in \cite{Belyakov2014}. Here we obtain approximations of both stability boundaries and study analytically both effects of stabilization and destabilization by damping, Moreover, the fourth approximation yields in addition boundaries of another stability domain. 

\section{Statement of the problem}
Consider linearization, $\dot{x}(t)=\mathbf J(t)\, x(t)$, of a nonlinear system about its equilibrium,
where $x(t)$ is the vector of state variable perturbations and $\mathbf J(t)$ is piecewise continuous, $T$-periodic and thus integrable \emph{Jacobian matrix} of the original nonlinear system. Solution of the matrix differential equation with the initial value being the identity matrix $\mathbf I$
\begin{equation}\label{eq:X}
\dot{\mathbf X}(t)=\mathbf J(t)\cdot\mathbf X(t), \quad \mathbf X(0)= \mathbf I,
\end{equation}
yields fundamental matrix and the monodromy matrix  $\mathbf F = \mathbf X(T)$.
If all eigenvalues of the monodromy matrix, \emph{Floquet multipliers}, have absolute values smaller than one, then the equilibrium of the nonlinear system is asymptotically stable, and if at least one eigenvalue has absolute value grater than one, then the equilibrium is unstable, see, e.g. \cite{SeyMai}. To  have analytical approximations of stability regions in parameter space we apply the following.  
\section{Averaging scheme}
Let the Jacobian matrix $\mathbf J(t)$ be expended into the series
\begin{equation}\label{eq:J}
\mathbf J(t) = \mathbf J_0(t) + \mathbf J_1(t) + \mathbf J_2(t)+\mathbf J_3(t)+\ldots,
\end{equation} 
where the lower index denotes the order of smallness. Suppose we know solution $\mathbf X_0(t)$ of the matrix initial value problem $\dot{\mathbf X}_0(t)=\mathbf J_0(t)\cdot\mathbf X_0(t)$, where $\mathbf X_0(0)= \mathbf I$.
Then the change of variable $\mathbf X(t) = \mathbf X_0(t) \cdot \mathbf Y(t)$ converts (\ref{eq:X}) to the standard form: 
\begin{equation}\label{eq:Y}
\dot{\mathbf Y}(t)=\mathbf H(t)\cdot\mathbf Y(t) , \quad \mathbf Y(0)= \mathbf I,
\end{equation}
where matrix $\mathbf H(t):=\mathbf X_0^{-1}(t)\cdot\left(\mathbf J(t)-\mathbf J_0(t)\right)\cdot \mathbf X_0(t)$ is small for $t\in[0,T]$.
Approximate solution of (\ref{eq:Y}) can be found with \emph{averaging method} as follows. 
Let 
\begin{equation}\label{eq:H}
\mathbf H(t) = \mathbf H_1(t) + \mathbf H_2(t)+\mathbf H_3(t)+\ldots,
\end{equation} 
where $\mathbf H_j(t):=\mathbf X_0^{-1}(t)\cdot \mathbf J_j(t)\cdot \mathbf X_0(t)$ for all $j\geq 1$. We will find solution as
\begin{equation}\label{eq:YE}
\mathbf Y(t) = \left(\mathbf I + \mathbf U_1(t) + \mathbf U_2(t)+\ldots \right)\cdot \mathbf Z(t),
\end{equation} 
where $\mathbf U_j(t)$ are matrix-functions, such that $\mathbf U_j(0) = \mathbf U_j(T)=0$ and $\mathbf Z(t)$ is the solution of the autonomous \emph{averaged  differential equation}:
\begin{equation}\label{eq:A}
\dot{\mathbf Z}(t)=\mathbf A \cdot\mathbf Z(t) , \quad \mathbf Z(0)= \mathbf I,
\end{equation}
where $\mathbf A = \mathbf A_1 + \mathbf A_2+\mathbf A_3+\ldots$, which can be written via the matrix exponential: %$\mathbf Z(t)=\exp\left(\mathbf A\, t \right)$:
\begin{equation}\label{eq:Z} 
\mathbf Z(t)=\exp\left(\left[\mathbf A_1 + \mathbf A_2+\mathbf A_3+\ldots \right] t \right).
\end{equation} 
The matrices $\mathbf A_j$ and matrix-functions  $\mathbf U_j(t)$
can be found one by one the following expressions.\footnote{The expressions are derived by 
differentiating (\ref{eq:YE}) w.r.t. time
\begin{align*}\dot{\mathbf Y}(t) = & \left(\dot{\mathbf U}_1(t) + \dot{\mathbf U}_2(t)+\ldots \right)\cdot \mathbf Z(t) + \left(\mathbf I + \mathbf U_1(t) + \mathbf U_2(t)+\ldots \right)\cdot \dot{\mathbf Z}(t),
\end{align*}
substituting there expressions for time derivatives from
(\ref{eq:Y}) and (\ref{eq:A}) 
\begin{align*}
 & \left(\mathbf H_1(t) + \mathbf H_2(t)+\ldots \right)\cdot \left(\mathbf I + \mathbf U_1(t) + \mathbf U_2(t)+\ldots \right)\cdot \mathbf Z  \\
& = \left(\dot{\mathbf U}_1(t) + \dot{\mathbf U}_2(t)+\ldots\right)\cdot\mathbf Z %\\ & 
+  \left(\mathbf I + \mathbf U_1(t) + \mathbf U_2(t)+\ldots \right) \cdot \left(\mathbf A_1 + \mathbf A_2+\ldots \right)\cdot\mathbf Z,
\end{align*}
collecting there terms of the
same order, and canceling non-degenerate matrix $\mathbf Z$, which yield the following equalities.
%\begin{enumerate} 

%\item[]
First order: 
$%\begin{align*}
\mathbf H_1(t)  = \dot{\mathbf U}_1(t) + \mathbf A_1. 
$% \end{align*}

%\item[]
Second order: 
$%\begin{align*}
\mathbf H_2(t) + \mathbf H_1(t)\cdot \mathbf U_1(t)  = \dot{\mathbf U}_2(t) + \mathbf U_1(t)\cdot\mathbf A_1 + \mathbf A_2, 
 $%\end{align*}

%\item[]
Third order: $ %\begin{align*} &
\mathbf H_3(t) + \mathbf H_1(t)\cdot \mathbf U_2(t)  + \mathbf H_2(t)\cdot \mathbf U_1(t)%\\ & 
  = \dot{\mathbf U}_3(t) + \mathbf U_2(t)\cdot\mathbf A_1 + \mathbf U_1(t)\cdot\mathbf A_2 + \mathbf A_3,
$ % \end{align*}

%\item[]
Fourth order:
$ %\begin{align*} &
\mathbf H_4(t) + \mathbf H_1(t)\cdot \mathbf U_3(t)  + \mathbf H_2(t)\cdot \mathbf U_2(t)+ \mathbf H_3(t)\cdot \mathbf U_1(t)
%\\  &
  = \dot{\mathbf U}_4(t) + \mathbf U_3(t)\cdot\mathbf A_1 + \mathbf U_2(t)\cdot\mathbf A_2  + \mathbf U_1(t)\cdot\mathbf A_3 + \mathbf A_4
$, %\end{align*}
%\item[]
and so on...
%\end{enumerate} 
}

For the first order approximation we calculate $\mathbf A_1$ as the average of $\mathbf H_1(t)$
\begin{equation}\label{eq_AVG1}
    \mathbf A_1 = \frac{1}{T}\int_0^{T} \mathbf H_1(t)  \dd t,%,
\end{equation}
under assumption that $\mathbf A_1$ and $\mathbf H_1(t)$ are of the same order of smallness. In particular we assume that $\mathbf H_1(t)$ does not contain periodic functions with small frequencies, which would appear in the denominator during integration and could cause high value of $\mathbf A_1$, thus violating the assumption of its smallness.

For the second order approximation we have to calculate
\begin{align}\label{eq_AVG2}
    \mathbf U_1(t) =  \int_0^t (&\mathbf H_1(\tau) - \mathbf A_1)  \dd \tau,\\ 
    \mathbf A_2 =  \frac{1}{T}\int_0^{T} (&\mathbf H_2(t) + \mathbf H_1(t)\cdot \mathbf U_1(t) -\mathbf U_1(t)\cdot\mathbf A_1)  \dd t,\label{eq_AVG2F}
\end{align}
using matrix $\mathbf A_1$ already obtained in (\ref{eq_AVG1}). 

For the third order approximation we have
\begin{align}
    \mathbf U_2(t) = \int_0^t (\mathbf H_2(\tau) - \mathbf A_2 %\nonumber\\[-10pt]
    +&\mathbf H_1(\tau)\cdot \mathbf U_1(\tau) - \mathbf U_1(\tau)\cdot \mathbf A_1) \dd \tau,\label{eq_AVG3}\\ 
     \mathbf A_3  =  \frac{1}{T}\int_0^{T} (\mathbf H_3(t) %\nonumber\\[-10pt] 
+ & \mathbf H_2(t)\cdot \mathbf U_1(t) - \mathbf U_1(t)\cdot\mathbf A_2 \nonumber\\[-8pt]
+ &\mathbf H_1(t)\cdot \mathbf U_2(t) - \mathbf U_2(t)\cdot\mathbf A_1)  \dd t.\label{eq_AVG3F}
\end{align}
%For the fourth order approximation we calculate
%\begin{align}
%\mathbf U_3(t) = \int\limits_0^t (&\mathbf H_3(\tau) - \mathbf A_3 \nonumber\\[-10pt] 
%+ & \mathbf H_2(\tau)\cdot \mathbf U_1(\tau) - \mathbf U_1(\tau)\cdot \mathbf A_2\nonumber\\ 
%+ & \mathbf H_1(\tau)\cdot \mathbf U_2(\tau) - \mathbf U_2(\tau)\cdot \mathbf A_1)   \dd \tau,\\
%\mathbf A_4  =  \frac{1}{T}\int\limits_0^{T} (\mathbf H_4(t) + &\mathbf H_3(t)\cdot \mathbf U_1(t)  - \mathbf U_1(t)\cdot\mathbf A_3\nonumber\\[-10pt] 
%+ &\mathbf H_2(t)\cdot \mathbf U_2(t)  - \mathbf U_2(t)\cdot\mathbf A_2\nonumber\\ +  & \mathbf H_1(t)\cdot \mathbf U_3(t)  -  %\mathbf U_3(t)\cdot\mathbf A_1) \dd t,\label{eq_AVG4F}
%\end{align}
and so on ...
for the $n+1$-th order approximation we calculate
\begin{align}
\mathbf U_{n}(t) = \int_0^t (\mathbf H_{n}(\tau) - \mathbf A_{n} + & \mathbf H_{n-1}(\tau)\cdot \mathbf U_1(\tau) - \mathbf U_1(\tau)\cdot \mathbf A_{n-1}\nonumber\\[-8pt]
\ldots &  \nonumber\\
+ & \mathbf H_1(\tau)\cdot \mathbf U_{n-1}(\tau) - \mathbf U_{n-1}(\tau)\cdot \mathbf A_1)   \dd \tau,\\
\mathbf A_{n+1}  =  \frac{1}{T}\int_0^{T} (\mathbf H_{n+1}(t) 
+ &\mathbf H_{n}(t)\cdot \mathbf U_1(t)  - \mathbf U_1(t)\cdot\mathbf A_{n}\nonumber\\[-8pt] 
% + &\mathbf H_{n-1}(t)\cdot \mathbf U_2(t)  - \mathbf U_2(t)\cdot\mathbf A_{n-1}\nonumber\\
\ldots  &  \nonumber\\  
 +  & \mathbf H_1(t)\cdot \mathbf U_{n}(t)  -  \mathbf U_{n}(t)\cdot\mathbf A_1) \dd t.\label{eq_AVGnF}
\end{align}

%\section{Lyapunov exponent approximation}
%Fundamental matrix of the linearized system at $t \to T$ can be expressed via product of matrix exponentials $$\mathbf X(T) = \exp(\mathbf J_0 T)\cdot \exp(\mathbf A[T]\, T),$$
%where we recall that $\mathbf A[T]$ may depend on $T$. The major \emph{Lyapunov exponent} 
%$$\lambda_1 = \limsup_{T\to\infty} \frac{1}{T}\ln ||\mathbf X(T)|| = \limsup_{T\to\infty} \frac{1}{T}\ln ||\exp(\mathbf J_0)^T\cdot \exp(\mathbf A[\infty])^T||,$$
%where we assume that limit $\mathbf A[\infty]$ exists. 
% and use  Lie product formula to cancel $T$.
% \emph{Spectrum of Lyapunov exponents} consists of eigenvalues of the matrix logarithm $$\frac{1}{2T}\ln(\exp(\mathbf A^{\prime}[\infty]) \cdot \exp(\mathbf J_0^{\prime})\cdot  \exp(\mathbf J_0)\cdot \exp(\mathbf A[\infty])),$$ where $^{\prime}$ denotes transposition. Thus if we take approximate $\mathbf A$, we can obtain approximate Lyapunov exponents.

\section{Monodromy matrix approximation}
Due to $\mathbf Y(T) = \mathbf Z(T)$ we have
\begin{align}
\mathbf F & =  \mathbf X(T) = \mathbf F_0\cdot \mathbf Y(T) %\nonumber\\ & 
=  \mathbf F_0 \cdot \mathbf Z(T) = \mathbf F_0\cdot \left(\mathbf I + \mathbf Z_1(T) + \mathbf Z_2(T)+\ldots\right),\label{eq:FG}
\end{align} 
where we denote $\mathbf F_0 := \mathbf X_0(T)$ as the zero order approximation of monodromy matrix, $\mathbf F \approx \mathbf F_0$. 
Hence, we can write
 expressions to find terms of the expansion $\mathbf F = \mathbf F_0 + \mathbf F_1 + \mathbf F_2 + \mathbf F_3 + \mathbf F_4 + \ldots$,
where $\mathbf F_j = \mathbf F_0 \cdot \mathbf Z_j(T)$ with $\mathbf Z_0(T) = \mathbf I$.  Expansion of the matrix exponential in (\ref{eq:Z}) yields expressions for $\mathbf F_j$ via $\mathbf A_k$, where $k\leq j$.

 For the first order approximation we have $\mathbf Z_1(T) = \mathbf A_1 T$, so that
\begin{equation}\label{eq_F1}
\mathbf F_1 = \mathbf F_0 \cdot \mathbf Z_1(T) = \mathbf F_0\cdot \mathbf A_1 T .
\end{equation}

For the second order approximation the
expansion of the matrix exponential in (\ref{eq:Z}) up to the second order terms yields  $\mathbf Z_2(T) =  \mathbf A_2 T + \frac{1}{2} \mathbf A_1^2 T^2$ and
\begin{equation}\label{eq_F2}
\mathbf F_2 = \mathbf F_0\cdot \left(\mathbf A_2 T + \frac{1}{2} \mathbf A_1^2 T^2\right) .
\end{equation}

For the third order approximation we have
\begin{align}
\mathbf F_3 =  \mathbf F_0\cdot \Big(\mathbf A_3 T+ \frac{1}{2}\left( \mathbf A_1\cdot \mathbf A_2 + \mathbf A_2\cdot \mathbf A_1\right) T^2 
 + \frac{1}{6} \mathbf A_1^3 T^3&\Big) ,\label{eq_F3}
\end{align}

For the fourth order approximation we calculate
\begin{align} 
\mathbf F_4  =  \mathbf F_0\cdot \Big(\mathbf A_4 T & + \frac{1}{2} \left(\mathbf A_1\mathbf A_3+\mathbf A_2^2+\mathbf A_3\mathbf A_1\right) T^2 \nonumber\\ &
+\frac{1}{6} (\mathbf A_1^2 \mathbf A_2 + \mathbf A_1\mathbf A_2 \mathbf A_1 + \mathbf A_2\mathbf A_1^2)T^3+\frac{1}{24} \mathbf A_1^4 T^4\Big),
\label{eq_F4}
\end{align}
and so on ...

\section{Stability conditions in 2-dimensional case}

The eigenvalues of the monodromy matrix, \emph{Floquet multipliers}, determine the stability of the solution $x=(0,0)'$ of the linearized system.
Since $\mathbf F$ is the $2\!\times 2$ matrix its eigenvalues can be found analytically as roots $\rho_{1}$ and $\rho_{2}$ of the characteristic polynomial:
\begin{equation}\label{eq:CP}
  \rho^2 - \tr(\mathbf F)\,\rho + \det(\mathbf F)=0.
\end{equation}
%\begin{displaymath}
%\rho_{1} = \frac{\tr(\mathbf F) + \sqrt{\tr(\mathbf F)^2 - 4 \det(\mathbf F)}}{2}, \quad
%\rho_{2} = \frac{\tr(\mathbf F) - \sqrt{\tr(\mathbf F)^2 - 4 \det(\mathbf F)}}{2}.
%\end{displaymath}
Stability conditions ($|\rho_{1}| \leq 1$ and $|\rho_{2}| \leq 1$) 
written in the case of real roots as $\rho\in [-1,1]$ and in the case of complex conjugate roots as $\rho_{1}\rho_{2} \leq 1$, with the use of (\ref{eq:CP}) and Vieta's formula $\rho_{1}\rho_{2} = \det(\mathbf F)$ correspondingly, take the form
\begin{equation}\label{eq:ASC}
    |\tr(\mathbf F)| -1 \leq \det(\mathbf F) \quad \text{and} \quad \det(\mathbf F) \leq 1,
\end{equation}
where for asymptotic stability all inequalities should be strict, see, e.g.~\cite{Kuznetsov}, p.~213. For instability, it is sufficient that %at least one of the eigenvalues has absolute value grater than one ($|\rho_{1}| > 1$ or $|\rho_{2}| > 1$), so that 
at least one of the conditions in (\ref{eq:ASC}) is violated.

Let us find approximations of stability conditions. Trace can be written as
\begin{equation}\label{eq:trF}
\tr\!\left(\mathbf F \right) =  \tr(\mathbf F_0) +
\tr(\mathbf F_1) + \tr(\mathbf F_2) + \tr(\mathbf F_3) + \tr(\mathbf F_4) + \ldots,
%\tr(\mathbf F_0\cdot \mathbf A_1)\, T + \tr(\mathbf F_0\cdot \mathbf A_2)\, T + \frac{1}{2} \tr(\mathbf F_0\cdot \mathbf A_1^2) \, T^2 + \ldots.
\end{equation}
with $\mathbf F_1$, $\mathbf F_2$, $\mathbf F_3$, and $\mathbf F_4$ calculated by  (\ref{eq_F1}), (\ref{eq_F2}), (\ref{eq_F3}), and (\ref{eq_F4}).

Notice that due to Liouville's formula, see, e.g. \cite{YakStarzh}, we have 
\begin{align*}
\det(\mathbf F) & = e^{\int^{T}_{0}\tr(\mathbf J(\tau)) \dd \tau} %\\ & 
= e^{\int^{T}_{0}\tr(\mathbf J_0(\tau)) \dd \tau} e^{\int^{T}_{0}\left(\tr(\mathbf J_1(\tau)) + \tr(\mathbf J_2(\tau))+\tr(\mathbf J_3(\tau))+\ldots\right) \dd \tau}.
\end{align*}
and  since $\mathbf F = \mathbf F_0\cdot \mathbf Z(T)$ and (\ref{eq:Z})  determinant can be written as
\begin{align*}
\det(\mathbf F) & = \det(\mathbf F_0) \det(\mathbf Z(T))%\\ & 
=  e^{\int^{T}_{0}\tr(\mathbf J_0(\tau)) \dd \tau} e^{\tr(\mathbf A_1)\, T + \tr(\mathbf A_2)\, T +\tr(\mathbf A_3)\, T+\ldots},
\end{align*}
so that for any dimension for all $j \geq 1$ we have
\begin{equation}\label{eq:trAj}
 \tr(\mathbf A_j) = \frac{1}{T}\int^{T}_{0}\tr(\mathbf J_j(\tau)) d \tau.
\end{equation}
Then expansion of the matrix exponential yields
\begin{align}
\det\!\left( \mathbf F \right) = e^{\int^{T}_{0}\tr(\mathbf J_0(\tau)) \dd \tau} \Big(1 & + \tr(\mathbf A_1) T + \tr(\mathbf A_2) T %\nonumber\\ & 
+ \frac{1}{2}(\tr(\mathbf A_1) T)^2 +\ldots\Big).\label{eq:detF}
\end{align}
The first order approximation of stability conditions $|\tr(\mathbf F)| -1 \leq \det(\mathbf F) \leq 1$ can be written with the use of (\ref{eq:trF}) and (\ref{eq:detF}) as
$$|\tr(\mathbf F_0) + \tr(\mathbf F_1)| -1 \,\leq\, e^{\int^{T}_{0}\tr(\mathbf J_0(\tau)) \dd \tau}\!\left(1 + \tr(\mathbf A_1) T\right)\, \leq\, 1$$
and so on ...

\section{Two DOF system with impulse parametric excitation}

\subsection{High frequency stabilization of inverted pendulum}

Let's show that motion of an inverted pendulum can be stable if we supply to the suspension point rather high frequency vibrations in vertical direction $l \ddot\varphi + b l \dot\varphi + \left(g \pm c\right) \sin(\varphi) = 0$, where $l$ is the length of the pendulum, $a\ll l$ is the amplitude of the vibrations of the pivot. The period of the pivot vibrations we normalize to be $T=2\pi$, besides in any semi-period acceleration of the pivot is constant and is equal $\pm c$, which sign changes each semi-period. We assume linear viscous friction with coefficient $b$.
It turns out that for rather low relative eigenfrequency $\omega  \ll 1 $ the inverted position becomes stable.   
The equation of motion can be written in the form 
\begin{align}
  \dot \varphi = & s,\label{eq:dphi1}\\
	\dot s = & - \beta\omega s - \left(\omega^2 \pm \eps\right) \sin(\varphi),\label{eq:ds1}
\end{align}
where $ \omega^2 = g/l$ is the relative eigenfrequency, $\eps=c/l$ is the relative excitation acceleration, and $\beta = b/g$ is the new damping coefficient.
Stability condition for this problem without damping can be found in \cite{ArnoldODE}. The linearized case without damping coincides with the Meissner equation, \cite{Meissner}.

 This system, linearized about inverted vertical position $(\phi,s)=(\pi,0)$, has the form $\dot{x}(t)=\mymat J(t)\, x(t)$, with vector $x(t)$ corresponding to the perturbation of vector $(\varphi(t),s(t))'$ and $\mymat J(t)$ being the piecewise constant Jacobian matrix of the original system: $\mymat J(t) = \mymat J^{+}$ if $t \in \left[0,\pi\right)$ and $\mymat J(t)=\mymat J^{-}$ if $t  \in \left[\pi,2\pi\right)$, where 
\begin{equation}\label{eq:Jpm}
\mymat J^{+}\! = \left(%
\begin{array}{c@{\quad}c}
  0 & 1 \\
  \omega^2  + \eps & - \beta\omega \\
\end{array}%
\right),\qquad \mymat J^{-}\! = \left(%
\begin{array}{c@{\quad}c}
  0 & 1 \\
  \omega^2  - \eps &  - \beta\omega\\
\end{array}%
\right).
\end{equation}
It is easy to find exact stability conditions to which approximate conditions converge as we shall demonstrate. 

\subsection{Exact stability conditions}
We have the following explicit expression of the monodromy matrix via matrix exponentials
\begin{equation}\label{eq:Fpm}
\mymat F = \exp\!\left(\pi\, \mymat J^{-}\right) \cdot \exp\!\left(\pi\, \mymat J^{+}\right).
\end{equation}

Since the determinant of a matrix product equals the product of the determinants, we have from (\ref{eq:Fpm})
\begin{equation}\label{eq:detFpm}
\det(\mymat F) = e^{\pi\tr(\mymat J^{-})} 
 e^{\pi\tr(\mymat J^{+})} = e^{-2\pi\beta\omega},
\end{equation} 
where we take into account that $\det(\exp\!\left(\pi\, \mymat J^{\pm}\right)) = e^{2\pi\tr(\mymat J^{\pm})}.$
The same expression as (\ref{eq:detFpm}), $\det(\mymat F) = e^{-2\pi\beta\omega }$, can be obtained by Liouville's formula for any piecewise continuous integrable $2\pi$-periodic modulation function, see \cite{YakStarzh}.
So with positive damping coefficient,
$\beta > 0$, asymptotic stability can only be lost when the first condition in 
(\ref{eq:ASC}) is violated, i.e. when
\begin{equation}\label{eq:ASC2}
    |\tr(\mymat F)| - 1 > e^{-2\pi\beta\omega}.
\end{equation} 
We compare exact stability borders, determined by (\ref{eq:Jpm}), (\ref{eq:Fpm}), and (\ref{eq:ASC2}) as
\begin{align}\label{eq:ASC3}
 \left|\tr\left(\exp \left[\pi%
\begin{pmatrix}%{c@{\quad}c}
   0 &  1 \\
   \omega^2\!-\!\eps &  - \beta\omega\\
\end{pmatrix}%
\right] \cdot \exp \left[\pi\left(%
\begin{matrix}%{c@{\quad}c}
  0 &  1 \\
  \omega^2 \!+\! \eps &  - \beta\omega \\
\end{matrix}%
\right)\right]\right)\right| - 1 = e^{-2\pi\beta\omega},
\end{align}
  with approximate stability boundaries obtained for this example.

\subsection{Approximate stability conditions}

We expend this matrix in the series $\mymat J^{\pm} = \mymat J_{0} + \mymat J_{1}^{\pm} +\mymat J_{2}$, where
\begin{equation}\label{eq:JKP01}
\mymat J_{0} = \!\left(%
\begin{array}{c@{\quad}c}
  0 & 1 \\
   0 & 0 \\
\end{array}%
\right)\!,\qquad \mymat J_{1}^{\pm} =\! \left(%
\begin{array}{c@{\quad}c}
  0 & 0 \\
   \pm \eps  & 0 \\
\end{array}%
\right)\!,\qquad \mymat J_{2} =\! \left(%
\begin{array}{c@{\quad}c}
  0 & 0 \\
   \omega^2  & - \beta\omega\\
\end{array}%
\right)\!,
\end{equation}
assuming that $\eps$, $\omega$, and $\beta$ have the same order of smallness.
\begin{equation}\label{eq:expJ0KP}
\mathbf X_0(t) = \exp(\mymat J_0 t) = \left(%
\begin{array}{c@{\quad}c}
  1 &  t \\
  0 & 1\\
\end{array}\right)\!,\quad
\mathbf X_0^{-1}(t) = \exp(-\mymat J_0 t) = \left(%
\begin{array}{c@{\quad}c}
  1 &  -t \\
  0 & 1 \\
\end{array}\right)\!.
\end{equation}
According the the formula $\mymat H_j(t):=\mathbf X_0^{-1}(t)\cdot \mymat J_j(t)\cdot \mathbf X_0(t)$ we have
\begin{equation}
 \mymat H_1^{\pm}(t)  =  \pm\eps \left(%
\begin{array}{c@{\quad}c}
  -t  &- t^2\\
  1 &  t \\
\end{array}\right)
\end{equation}
Formula (\ref{eq_AVG1}) reads as
\begin{align}
    \mymat A_1 = \frac{1}{2\pi}\int\limits_0^{\pi} \mymat H_1^{+}(t) \dd t + \frac{1}{2\pi}\int\limits_{\pi}^{2\pi} \mymat H_1^{-}(t) \dd t %\nonumber\\ 
    = \frac{\eps}{2\pi} \left(%
\begin{array}{c@{\quad}c}
  \pi^2 &  2\pi^3\\
  0 & -\pi^2 \\
\end{array}\right),\label{eq_KP_AVG1}
\end{align}
The zero order approximation of monodromy matrix is the following
\begin{equation}\label{eq:expJ0KP_T}
\mymat F_0 = \mymat X_0(2\pi) =\exp(\mymat J_0 2\pi) = \left(%
\begin{array}{c@{\quad}c}
  1 &  2\pi \\
  0 & 1 \\
\end{array}\right).
\end{equation}
With (\ref{eq_F1}) we find the first order adjustment of monodromy matrix 
\begin{equation}
\mymat F_1 = \mymat F_0\cdot \mymat A_1 2\pi =  \pi^2\eps\left(%
\begin{array}{c@{\quad}c}
  1 &  0\\
  0 & -1 \\
\end{array}\right).
\end{equation}
Thus $\tr (\mymat F_0) = 2$, $\tr (\mymat F_1) = 0$, and $\tr (\mymat A_1) = 0$.

For the second order approximation we take
$%\begin{equation}
 \mymat H_2(t)  =  \left(%
\begin{array}{c@{\quad}c}
  -t \omega^2  &- t\left(\omega^2 t - \beta\omega\right)\\
  \omega^2 &  \omega^2 t - \beta\omega \\
\end{array}\right)
$ %\end{equation}
and obtain with (\ref{eq_AVG2})
%$$
%\mymat U_{1}^{\pm} = \eps\left(%
%\begin{array}{c@{\quad}c}
%  -\frac{\left(\pi\pm\left( t - \pi\right)\right)\left(t+\pi\right)}{2} & -\pi^2 t - \frac{\pi^3 \pm (t^3 - \pi^3)}{3} \\
%   \pi \pm \left( t - \pi\right)  & \frac{\left(\pi\pm\left( t - \pi\right)\right)\left(t+\pi\right)}{2} \\
%\end{array}%
%\right).
%$$
and (\ref{eq_AVG2F}) the matrix
$$ \mymat A_2 = \frac{1}{2\pi} \left[ \begin {array}{c@{\quad}c} 
 \frac{2}{3}\,{\eps}^{2}{\pi}^{4}-2{\pi}^{2}{\omega}^{2}&{
\frac {4\,{\eps}^{2}{\pi}^{5}}{15}}-\frac{8}{3}\,{\pi}^{3}{\omega}^{2}+2{\pi}^
{2}\,\beta\,\omega\\[5pt] 
-\frac{2}{3}\,{\eps}^{2}{\pi}^{3}+2\pi\,{\omega}^{2}&-\frac{2}{3}\,{
\eps}^{2}{\pi}^{4}+2{\pi}^{2}{\omega}^{2}-2\pi\,\beta\,\omega\end {array} \right]
.
$$
We calculate according (\ref{eq_F2}) the second order adjustment of monodromy matrix 
\begin{equation}
    \mymat F_2 = \left[ \begin {array}{c@{\quad}c} -\frac{1}{6}\,{\pi}^{4}{\eps}^{2}+2{
\pi}^{2}{\omega}^{2}& -\frac{1}{15}\,{\pi}^{5}{\eps}^{2}+\frac{4}{3}\,{\pi}^{3}{\omega}^{2}-2 {\pi}^{2}\beta\,\omega\\[5pt]
-\frac{2}{3}\,{\pi}^{3}{\eps}^{2}+2\pi\,{\omega}^{2} & -\frac{1}{6}\,{\pi}^{4}{\eps}^{2}+2{\pi}^{2}{\omega}^{2}-2\pi\,\beta
\,\omega\end {array} \right],
\end{equation}
Thus $\tr (\mymat F_2) = -\frac{1}{3}\,{\pi}^{4}{\eps}^{2}+4{\pi}^{2}{\omega}^{2}-2\pi\,\beta
\,\omega$ and $\tr (\mymat A_2) = -2\pi\,\beta\,\omega$.

Second order approximation of stability border written from (\ref{eq:ASC2}) as $$|\tr(\mymat F_0)+\tr(\mymat F_1)+\tr(\mymat F_2)| - 1 = 1 - 2\pi\beta\omega$$ yields $4\,{\pi}^{2}{\omega}^{2}-\eps_p^{2}{\pi}^{4}/3=0$ and $4-4\pi\,\beta\,w+4{\pi}^{2}{\omega}^{2}-\eps_n^{2}{\pi}^{4}/3=0$ in cases of positive and negative value of the sum $\tr(\mymat F_0)+\tr(\mymat F_1)+\tr(\mymat F_2)$ correspondingly. Hence we have corresponding stability borders denoted by indexes $p$ and $n$
$$\eps_p = \frac{2\sqrt{3}}{\pi}\,\omega,\quad \eps_n = \frac{2\sqrt{3}}{\pi}\sqrt{\omega^2 - \frac{\beta\omega}{\pi}+ \frac{1}{\pi^2}}, $$
same as in the third approximation, see dashed lines in the Figure on the left, %~\ref{f:InvPendStab}.
because  (\ref{eq_AVG3})--(\ref{eq_F3}) yield $\tr (\mymat F_3) = 0$ and $\tr (\mymat A_3) = 0$.
 
Fourth order approximation, where stability boundary equation reads as 
\begin{align*}
|\tr(\mymat F_0)+\tr(\mymat F_1)+\tr(\mymat F_2)+\tr(\mymat F_3)+\tr(\mymat F_4)| - 1 %\nonumber\\ 
= 1 - 2\pi\beta\omega + \frac{1}{2} \left(2\pi\beta\omega\right)^2,
\end{align*}
  yields in cases of positive and negative sums of traces the following two equations:\\
$%\begin{align*}
\footnotesize\frac{{\pi}^{8}{\eps}^{4}_p}{1260}- \frac{{\pi}^{4}}{3}\!\left( 1+{\frac {4{\pi}^{2}{\omega}^{2}
}{15}}-{\pi} \beta \omega \right) {\eps}^{2}_p %\\
+4{\pi}^{2}{\omega}^{2}\left(1+\frac{{\pi}^{2
}{\omega}^{2}}{3}-\beta{\omega}{\pi}\right)=0,
$\\%\end{align*} 
$%\begin{align*}
\footnotesize\frac {{\pi}^{8}{\eps}^{4}_n}{1260}- \frac{{\pi}^{4}}{3}\!\left( 1+{\frac {4{\pi}^{2}{\omega}^{2}
}{15}}-{\pi} \beta \omega \right) {\eps}^{2}_n%\\
+4\left(1-\beta \pi \omega+{\pi}^{2}{\omega}^{2}\!\left(1+\frac{{\pi}^{2
}{\omega}^{2}}{3}-\beta{\omega}{\pi}+{\beta}^{2}\right)\right)\!=\!0.
$\\ %\end{align*} 
Solutions of these two equations with respect to $\eps^2$ give us four borders, drawn in the Figure with solid lines,  which approximate two exact stability domains determined by (\ref{eq:ASC3}) and marked in gray.

\section{Conclusion}

We develop convenient algorithm for obtaining approximate stability boundaries of parametrically excited systems.
We demonstrate how this algorithm can be applied to the particular case of parametric pendulum for obtaining approximate stability  domains even in the case of damping and impulse parametric excitation. Stabilizing and destabilizing effects of damping on inverted equilibrium of the parametric pendulum are revealed with the use of the fourth order approximation of stability boundaries.
    
\begin{figure}
\centering
\includegraphics[width=0.49\columnwidth]{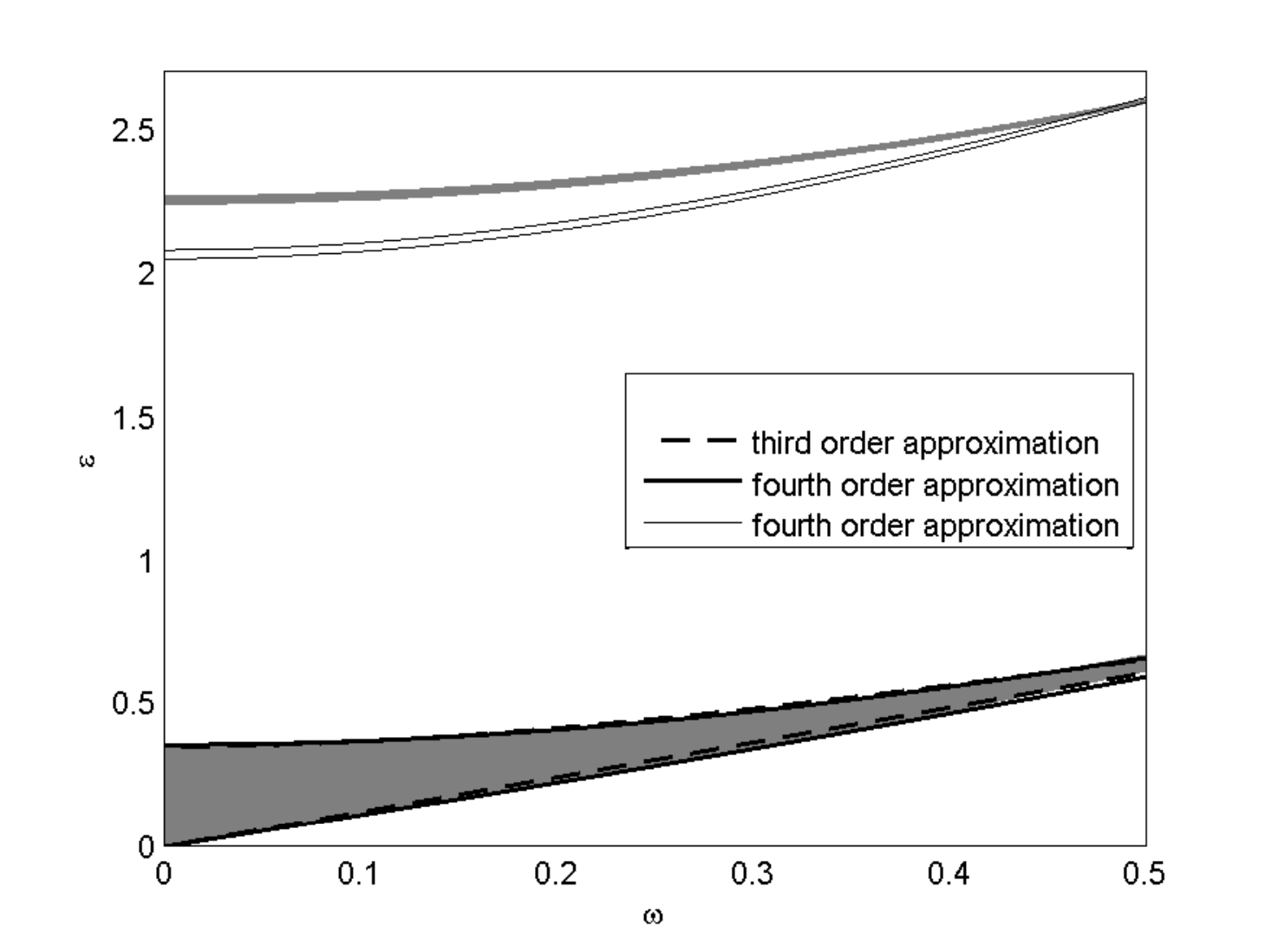}
\includegraphics[width=0.49\columnwidth]{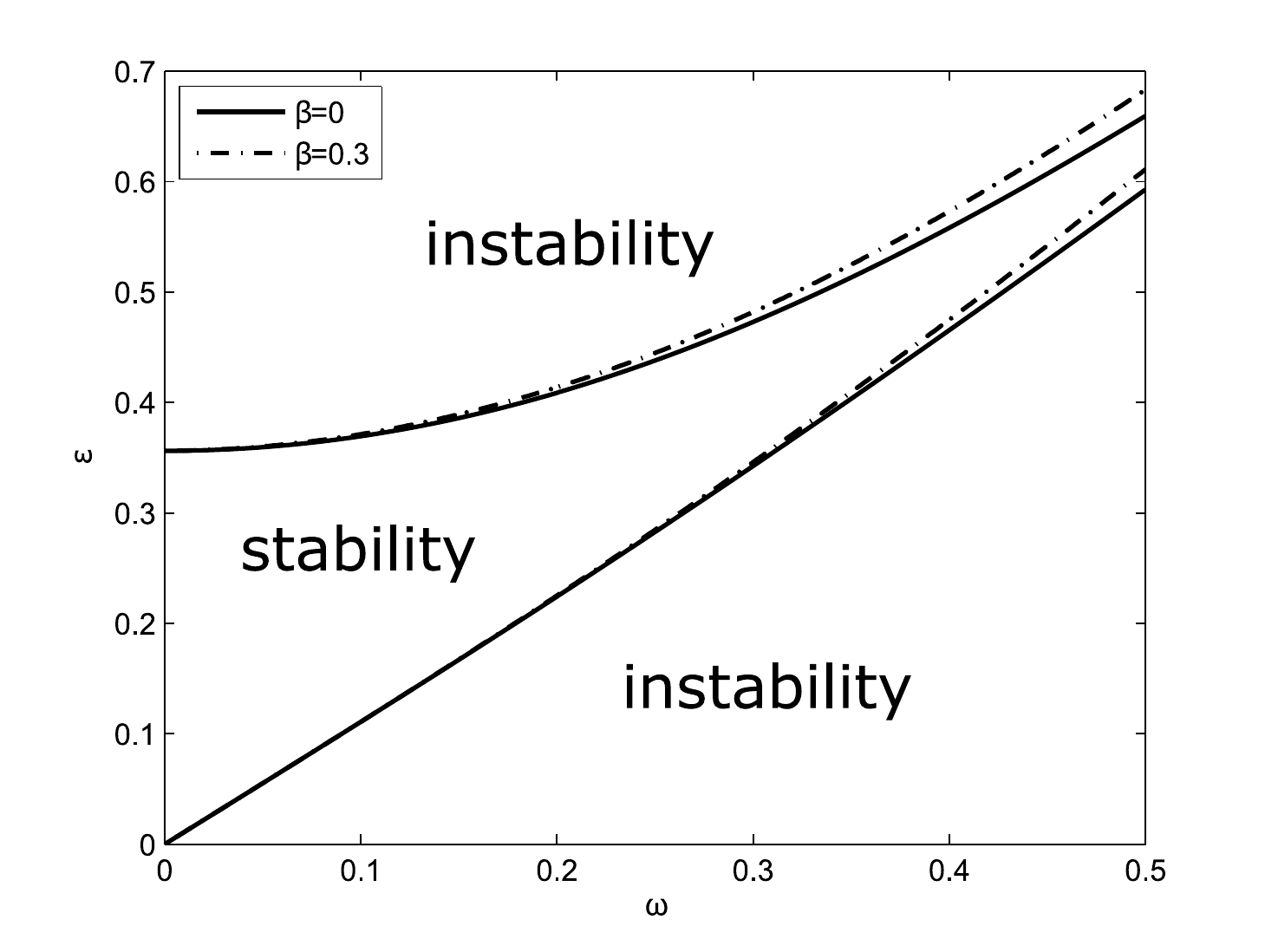}
\caption{Comparison of approximate stability boundaries of the third (dashed lines) and fourth (solid lines) orders with exact stability domains (gray areas) on the left. Damping stabilization and destabilization of the inverted vertical pendulum position on the right.}\label{f:Q}
\end{figure}
In the right figure we draw approximate stability boundaries of inverted vertical pendulum position. Addition of small linear viscous friction $\beta$ shifts both stability boundaries upward. Thus, at the lower boundary additional friction destabilizes the inverted pendulum while at the upper boundary friction stabilizes the pendulum position.
% Detailed analysis will be given in the extended paper.

\section{Acknowledgments}
%A.P. Seyranian received funding from the Russian Foundation for Basic Research grant 19-01-00660.\\
%A.O. Belyakov received funding from the Russian Foundation for Basic Research grant 18-010-01169.
A.O. Belyakov received funding from the Russian Science Foundation grant 19-11-00223.


\end{document}